\documentclass[12pt,oneside]{amsart}
\usepackage[a4paper,nofoot]{geometry}
\let\classno\subjclass
\renewcommand\bibname{\def\and{{\normalfont\rmfamily and }}%
                    \normalfont\scshape}
\usepackage{amsmath,amssymb} 
\newtheorem{thm}{Theorem}         
\newtheorem{prop}[thm]{Proposition} 
\newtheorem{lem}[thm]{Lemma}

\let\nc\newcommand  \let\dmo\DeclareMathOperator
\makeatletter
\def\bb@symb#1\f@@{\expandafter\nc\csname bb#1\endcsname{\mathbb{#1}}}
\nc\bbsymbols[1]{
  \@for\@current:=#1\do{\expandafter\bb@symb\@current\f@@}}
\def\cal@symb#1\f@@{\expandafter\nc\csname c#1\endcsname{\mathcal{#1}}}
\nc\calsymbols[1]{
  \@for\@current:=#1\do{\expandafter\cal@symb\@current\f@@}}
\def\frak@symb#1\f@@{\expandafter\nc\csname frak#1\endcsname{\mathfrak{#1}}}
\nc\fraksymbols[1]{
  \@for\@current:=#1\do{\expandafter\frak@symb\@current\f@@}}
\def\mym@th@p#1\f@@{\expandafter\dmo\csname#1\endcsname{#1}}
\nc\operators[1]{
  \@for\@current:=#1\do{\expandafter\mym@th@p\@current\f@@}}
\makeatother

\bbsymbols{Z,Q,C}
\let\Z\bbZ\let\Q\bbQ
\fraksymbols{sl,g,a,sp,so}
\let\g\frakg\let\sl\fraksl
\nc\ax{\fraka_x}
\operators{Gal,Lie,Sym,ad,Aut}
\nc\G{\Gamma}
\nc\Qbar{\overline{\Q}}
\nc\Ql{\Q_\ell}
\nc\prim{{\mathrm{prim}}}
\nc\Frob{\mathrm{Frob}}
\nc\abs[1]{\left|#1\right|}

\title[$\ell$-adic representations attached to noncongruence subgroups]%
 {On some $\ell$-adic representations of $\Gal(\Qbar/\Q)$
   attached to noncongruence subgroups.}  

\author{A. J. Scholl}

\classno{11F80 (Primary); 11G18 (Secondary)}

\begin{document}
\begin{abstract}
  The $\ell$-adic parabolic cohomology groups attached to
  noncongruence subgroups of $SL_2(\Z)$ are finite-dimensional
  $\ell$-adic representations of $\Gal(\Qbar/K)$ for some number field
  $K$. We exhibit examples (with $K=\Q$) for which the primitive parts
  give Galois representations whose images are open subgroups of the
  full group of symplectic similitudes (of arbitrary dimension). The
  determination of the image of the Galois group relies on Katz's
  classification theorem for semisimple subalgebras of $\sl_n$
  containing a principal nilpotent element, for which we give a short
  conceptual proof, suggested by I. Grojnowski.
\end{abstract}
\maketitle
\thispagestyle{empty}

\section{Introduction}

Let $\G\subset PSL_2(\Z)$ be a subgroup of finite index. In the papers
\cite{S1,S2,S3} we studied $\ell$-adic Galois
representations attached to cusp forms on $\G$. Attached to $\G$ is
a certain field $K_\G$ and, for each even integer $k\ge0$, a
compatible system of $\ell$-adic representations
\[
\rho_\ell=\rho_{\ell,k,\G}\colon\Gal(\Qbar/K_\G) \to GL_{2d}(\Ql)
\]
where $d=d_{k+2}$ is the dimension of the space $S_{k+2}(\G)$ of
cusp forms on $\G$ of weight $k+2$. These representations are
defined using $\ell$-adic parabolic cohomology, and are a mild
generalisation of the $\ell$-adic representations of Deligne
\cite{De}. If $\G'$ is the smallest congruence subgroup of $SL_2(\Z)$
containing $\G$ then $\rho_{\ell,k,\G}$ contains as an invariant
subspace the restriction of $\rho_{\ell,k,\G'}$ to
$\Gal(\Qbar/K_\G)$. The representation we are concerned with here
is the quotient, which we denote $\rho_{\ell,k,\G}^\prim$.

For any $\G$ the representations $\rho_{\ell,k,\G}$ are the
$\ell$-adic realisations of a certain motive (in the sense of
Grothendieck) $M_{k,\G}$ defined over $K_\G$. (For congruence
subgroups this was shown in \cite{S2.5}, and the trivial
generalisation to other groups was explained in \cite{S3}.) The
Hodge type of $M_{k,\G}$ is of the form $(k+1,0)^d+(0,k+1)^d$,
and so the representations $\rho_{\ell,k,\G}$ are (by Faltings
\cite{Fa}) Hodge-Tate of the same type. Moreover by Deligne's
proof of the Weil conjectures, they are pure of weight $k+1$. As
a final general remark, there is a perfect pairing
\[
M_{k,\G}\otimes M_{k,\G} \to \Q(-k-1)
\]
which is alternating (since $k$ is even) and so the image of
$\rho_{\ell,k,\G}$ is (after suitable conjugation) contained in
$GSp_{2d}(\Ql)$, the group of symplectic similitudes. The same
statements hold for the quotient $\rho_{\ell,k,\G}^\prim$ (since it
is the kernel of an algebraic projector, given by the
trace from $\G'$ to $\G$).

We considered in \cite{S3} the following three subgroups of
$PSL_2(\Z)$.  Write $\G_{43}$ for the subgroup generated by the
matrices
\[
\begin{pmatrix}1&4\\0&1\end{pmatrix},\quad
\begin{pmatrix}2&1\\1&1\end{pmatrix},\quad
\begin{pmatrix}1&-1\\2&-1\end{pmatrix}
\]
and $\G_{52}$ for the subgroup generated by
\[
\begin{pmatrix}1&5\\0&1\end{pmatrix},\quad
\begin{pmatrix}0&-1\\1&0\end{pmatrix},\quad
\begin{pmatrix}2&3\\1&2\end{pmatrix}.
\]
Thus $\G_{43}$ and $\G_{52}$ both have index 7 and two cusps, of
widths 4 and 3 (5 and 2, respectively). Also let $\G_{711}$ be
the subgroup of index 9 generated by
\[
\begin{pmatrix}1&7\\0&1\end{pmatrix},\quad
\begin{pmatrix}0&-1\\1&0\end{pmatrix},\quad
\begin{pmatrix}3&-4\\1&-1\end{pmatrix},\quad
\begin{pmatrix}-1&-4\\1&3\end{pmatrix}
\]
which has a cusp of width 7 and two cusps of width 1. 

If $\G$ is one of these three groups then it can be shown (cf.\
\cite{S3},  \S4.9) that $K_\G=\Q$. By applying standard formulae
for the dimensions of spaces of modular forms, we find that in
each case $\dim\rho_{\ell,k,\G}^\prim=k$.

Using methods from algebraic geometry, and in particular the theory of
vanishing cycles, we obtained in \cite{S3} a criterion for the image
of $\rho_\ell$ to contain a unipotent element with a ``long'' Jordan
block.  In particular, in \S4 of \cite{S3} the following result is
proved:

\begin{thm}
Let $\G$ be one of $\G_{52}$, $\G_{43}$, $\G_{711}$. Let $p=7$,
$7$ or $2$ respectively, and let $\ell\ne p$.  Let $k\ge2$ be even.
Then the image under $\rho_{\ell,k,\G}^\prim$ of an inertia subgroup at $p$
contains a unipotent element $X$ such that $(X-1)^{k-1}\ne0$.  
\end{thm}

We now fix once and for all a prime $\ell$ different from the prime
$p$ of Theorem~1, and write $C$ for the completion of the algebraic
closure of $\Ql$.  Let $G_{k,\G}\subset GSp_{k/C}$ be the connected
component of the identity in the Zariski closure of the image of
$\rho_{\ell,k,\G}^\prim$. It is a connected algebraic group over $C$.
In this paper we use Theorem 1 to prove:

\begin{thm}
Let $\G$ be as in Theorem~1, and $k\ge2$ an even integer. Then
$G_{k,\G}= GSp_{k/C}$.
\end{thm}

By Bogomolov's theorem \cite{Bog} it follows that the image of
$\rho_{k,\G}^\prim$ is an open subgroup of $GSp_k(\Q_\ell)$.

Apart from showing that the motives associated to non-congruence
subgroups can in some sense be as general as possible, Theorem 2 also
gives an explicit construction, for every even $k$ and every prime
$\ell$, of an $\ell$-adic representation $\rho\colon\Gal(\Qbar/\Q)\to
GSp_k(\Ql)$ with open image, which occurs in the $\ell$-adic
cohomology of a smooth projective variety over $\Q$. It does not seem
easy to produce examples of such representations by other methods.

These methods apply also to the case of $k$ odd (although there is
some ambiguity in the notion of field of definition for odd weight ---
see \cite[Remark 5.10(iii)]{S1} for a discussion) and, although we
have not checked all the details, it seems likely that one will obtain
odd-dimensional representations of $\Gal(\Qbar/\Q)$ whose image is
open in a group of orthogonal similitudes (except perhaps in the case
$k=7$, where a group of type $G_2$ might conceivably occur).

%% Summary of Bogomolov, CRAS 280 (1980) 701:
%%
%% $A/K$ abelian variety over number field. The image of Galois in the
%% $\ell$-adic representation is an open subgroup of an algebraic
%% subgroup of $GL_{2g}$, its Lie algebra is algebraic and contains the
%% homotheties. (Corollary is weak Manin-Mumford.)
%%
%% Suppose $\rho\colon G_K\to \Aut(V)$ is a (global) $\ell$-adic
%% representation. Let $G_V$ be its image, $\g_V=\Lie G_V$. Also $H_V=$
%% algebraic envelope of $G_V$, $\h_V=\Lie H_V$. Suppose $V$ is
%% Hodge-Tate at every place dividing $\ell$. Then $G_V$ is open in
%% $H_V$, and so $\g_V=\h_V$.
%%
%% Useful references: Serre in `\emph{Kyoto Sumposium on ALgebraic Number
%% Theory} 1977; also Serre in Ast\'erisque 65, 1979

\section{Number-theoretic part}

In this section we reduce Theorem~2 to a Lie-theoretic statement.  It
is convenient to axiomatise the properties of $\rho_{\ell,k,\G}^\prim$ we use.
Assume that we have a $\Ql$-vector space $V$ of dimension $k\ge2$ and
a continuous representation $\rho\colon \Gal(\Qbar/\Q) \to \Aut(V)$.
Let $G\subset GL(V)$ be the connected component of the identity in the
Zariski closure of the image of $\rho$.  Consider the following
conditions on $(\rho,V)$:

\begin{itemize}
\item[(H1)]  $\rho$ is pure of some weight $w\in\Z$;
\item[(H2)] The restriction of $\rho$ to $\Gal(\Qbar_\ell/\Ql)$ is
  Hodge-Tate, with exactly two Hodge-Tate weights;
\item[(H3)] For some $p\ne\ell$ there is an open subgroup
  $I'\subset I_p$ of the inertia group at $p$ such that the
  restriction of $\rho$ to $I'$ is unipotent and indecomposable.
\end{itemize}

We remind the reader: (H1) means that $\rho$ is
unramified outside a finite set $S$ of primes, and that for all
$p\notin S\cup\{\ell\}$ the eigenvalues of a geometric Frobenius at
$p$ are algebraic numbers, all of whose conjugates have absolute value
$p^{w/2}$. As for (H2), write $\sigma$ for the unique continuous
action of $\Gal(\Qbar_\ell/\Ql)$ on $C$ extending the Galois action on
$\Qbar_\ell$. Let $\chi_{\mathrm{cycl}}\colon
\Gal(\Qbar/\Q)\to\Z_\ell^*$ be the cyclotomic character, so that for
any $\ell^n$-th root of unity $\eta\in\Qbar$ and $g\in\Gal(\Qbar/\Q)$,
$g(\eta)=\eta^{\chi_{\mathrm{cycl}}(g)}$, and set
\[
V(i)=\bigl\{ v\in V\otimes_{\Ql} C \bigm|
(\rho\otimes\sigma)(g)v=\chi_{\mathrm{cycl}}^i(g)v \text{ for all
}g\in\Gal(\Qbar_\ell/\Ql) \bigr\}.
\]
Then $V$ is Hodge-Tate if the natural map $\bigoplus V(i)\otimes_{\Ql} C
\to V\otimes_{\Ql} C$ is an isomorphism, and its Hodge-Tate weights are those
$i$ for which $V(i)\ne0$. Finally, by Grothendieck's $\ell$-adic
monodromy theorem, (H3) is equivalent to the existence of some
$X\in \rho(I_p)$ whose Jordan form has a single block.

In the case $\rho=\rho_{k,\G}^\prim$, as explained in the
Introduction, (H1) is satisfied with $w=k+1$ and (H2) with weights
$\{0,-k-1\}$. Condition (H3) is the content of Theorem~1.

\begin{prop}
  Let $\rho\colon \Gal(\Qbar/\Q)\to \Aut(V)$ be a representation
  satisfying (H1) and (H3), and whose restriction to
  $\Gal(\Qbar_\ell/\Ql)$ is Hodge-Tate.  Then the restriction of
  $\rho$ to any open subgroup of $\Gal(\Qbar/\Q)$ is absolutely
  irreducible.
\end{prop}

\begin{proof}
  Let $E/\Q_p$ be a finite extension such that $I'$ is the inertia
  subgroup of $\Gal(\Qbar_p/E)$, let $q$ be the order of the residue
  field of $E$, and let $\Frob_q\in\Gal(\Qbar_p/E)$ be a geometric
  Frobenius --- that is, the inverse of any element lifting the
  $q$-power Frobenius on the residue field. In particular, 
  $\chi_{\mathrm{cycl}}(\Frob_q)=q^{-1}$.
  
  Hypothesis (iii) says that the Jordan normal form of
  $\rho\vert_{I'}$ has one block, so the invariants $V^{I'}$ form a
  1-dimensional subspace of $V$, on which $\Frob_q$ acts as a scalar
  $\alpha\in \Ql^*$. By the structure of the tame inertia group, the
  complete set of eigenvalues of $\rho(\Frob_q)$ is therefore
  $\{\alpha q^{j}\mid 0\le j\le k-1\}$.
  
  Recall that if $\chi\colon \Gal(\Qbar/\Q)\to \Ql^*$ is a continuous
  homomorphism whose restriction to the decomposition group at $\ell$
  is Hodge-Tate, then $\chi$ is the product of an integral power of
  $\chi_{\mathrm{cycl}}$ and a character of finite order. So there
  exists an integer $m$ and a character $\epsilon$ of finite order
  such that $\det\rho=\chi_{\mathrm{cycl}}^{-m}\epsilon$. By
  hypothesis (i), $m=wk/2$. Then
  \[
  \det\rho(\Frob_q)=\prod_{j=0}^{k-1}\alpha q^j=q^m\epsilon(\Frob_q)
  \]
  and so $\alpha$ is the product of $q^{(w-k+1)/2}$ and a root of
  unity.

  Now if %$\rho'$ 
  $V'\subset V$ is a $\Gal(\Qbar/\Q)$-invariant subspace
  %%subrepresentation of $\rho$ 
  of dimension $k'>0$
  then $V$ and $V'$ have the same space of $I'$-invariants
  %%we have $V'{}^{I'}=V^{I'}$ %$\rho'{}^{I'}=\rho{}^{I'}$ 
  (since $V^{I'}$ is $1$-dimensional) and $V'$ satisfies the
  hypotheses of the Proposition. Therefore the previous argument applied
  to $V'$ gives
  \[
  \abs{\alpha}=q^{(w-k'+1)/2},\qquad\emph{i.e.}\quad k'=k.
  \]
  So $\rho$ is irreducible. Finally, let $U$ be any subspace of the
  space of $\rho$ which is invariant under some open subgroup
  $H\subset\Gal(\Qbar/\Q)$, and let $g\in\Gal(\Qbar/\Q)$. Then $U$ and
  $\rho(g)U$ are both invariant under the open subgroup $I''=H\cap
  H^g\cap I'$ of $I_p$. But since the action of $I''$ also has one
  Jordan block, it has a unique invariant subspace of each dimension.
  So $\rho(g)U=U$, hence $U$ is invariant under $\Gal(\Qbar/\Q)$. So
  as $\rho$ is irreducible, its restriction to $H$ is also
  irreducible.
  
  Finally, the same argument carries through if we replace $\Ql$  by a
  finite extension, so the restriction of $\rho$ to any open subgroup
  is absolutely irreducible.
\end{proof}

As a consequence, since $G$ contains the image of an open subgroup
of $\Gal(\Qbar/\Q)$, it acts (absolutely) irreducibly on $V$, and
therefore (being connected by definition) it is reductive.  In
particular, for $k=\dim(V)=2$ we have $G=GL_2$. Henceforth we assume
that $(\rho,V)$ satisfies hypotheses (H1)--(H3) above, and that $k>2$.

Let $\g= \Lie G \cap \sl(V\otimes C)$; since $V$ is an irreducible
$G$-module, the centre of $\Lie G$ has dimension at most one, hence
$\g$ is a semisimple Lie algebra over $C$ .  By (H3) there exists a
unipotent element of $G$ whose Jordan decomposition has one block; let
$x\in\g$ be its logarithm.  Then $x$ is a nilpotent element of $\g$
which has just one Jordan block, viewed as an endomorphism of
$V\otimes C$.

Now recall the 1-dimensional Hodge-Tate torus associated to $\rho$ (as
a representation of the local Galois group).  Let $H_\ell\subset
GL(V\otimes_{\Ql} C)$ be the Zariski closure of the the image of
$\Gal(\Qbar_\ell/\Ql)$ by $\rho$. Since $\rho$ is Hodge-Tate, there is
a unique homomorphism $\zeta\colon \mathbb{G}_m \to H_\ell$ for which
$V(i)\otimes C$ is the eigenspace of the character $t\mapsto t^i$ of
$\mathbb{G}_m$. (See \cite{Se}, \S1.4, where $\zeta$ is denoted
$h_V$.) Passing to the Lie algebra, there is a unique semisimple
element $z_{HT}=d\zeta\in \Lie H_\ell\subset \Lie G$ such that
$V(i)\otimes C=\ker(z_{HT}-i)\subset V\otimes C$.

%% By the Jacobson-Morosov theorem, there is a
%% homomorphism $\lambda_{x}\colon \sl_2 \to\g$ such
%% \[
%% \lambda_{x}\colon\begin{pmatrix}0&1\\0&0\end{pmatrix}\mapsto x.
%% \]
%% Write $\ax\subset\g$ for the image of $\lambda_{x}$, and
%% \[
%% y=\lambda_{x}\begin{pmatrix}0&0\\-1&0\end{pmatrix},\quad
%% h=\lambda_{x}\begin{pmatrix}1&0\\0&-1\end{pmatrix}.
%% \]
%% (In the terminology of \cite{Bou}, VIII.11.1, $(x,h,y)$ is an
%% $\sl_2$-triplet.) 
%% 
%% Since $x$ has one Jordan block, the restriction of $\rho$ to $\ax$ is
%% isomorphic to the irreducible representation $\Sym^{k-1}$ of $\sl_2$,
%% and so $h$ has eigenvalues $\{k-1,k-3,\dots,-k+1\}$. Therefore any
%% non-zero semisimple element of $\ax$ has $k$ eigenvalues acting on $V$
%% (being a conjugate of non-zero multiple of $h$).
%% 
%% On the other hand, by hypothesis (ii) the Hodge-Tate element
%% $z_{HT}\in\Lie G$ has exactly 2 eigenvalues, namely the (integral)
%% Hodge-Tate weights of $\rho$. Therefore for some $a$, $b\in\Z$,
%% $az_{HT}+b$ is a semisimple element of $\g$ with exactly 2
%% eigenvalues. Therefore for $k>2$, $\g$ strictly contains $\ax$, and so
%% is not isomorphic to $\sl_2$.
 
We now appeal to the following result of Katz (the Classification
Theorem 9.10 in \cite{Ka}):

\begin{prop}
  Let $V$ be a finite-dimensional vector space of dimension $k$ over
  an algebraically closed field of characteristic zero, and $\g$ a
  semisimple Lie subalgebra of $\sl(V)$. Assume that $\g$ contains a
  nilpotent element $x$ which as an endomorphism of $V$ has only one
  Jordan block.  Then one of the following holds:
\begin{itemize}
\item[(i)] $\g\simeq\sl_2$ with $V\simeq \Sym^{k-1}$.

\item[(ii)] $\g=\sl(V)$;

\item[(iii)] $k$ is even, $\g=\fraksp(V)$ for a nondegenerate
alternating form on $V$;

\item[(iv)] $k$ is odd, $\g=\frakso(V)$ for a nondegenerate
symmetric form on $V$;

\item[(v)] $k=7$ and $\g$ is of type $G_2$.
\end{itemize}
\end{prop}

%% We have just excluded (assuming $k>2$) the case (i). In the setting of
%% Theorem~2, $k$ is even and $V$ has a $\g$-invariant symplectic form,
%% so we must be in case (iii), proving Theorem~2.

The hypothesis (H2) enables us to eliminate the case (i) unless
$k=2$. Indeed, the Hodge-Tate element
$z_{HT}\in\Lie G$ then has exactly 2 eigenvalues, namely the (integral)
Hodge-Tate weights of $\rho$. Therefore for some $a$, $b\in\Z$,
$az_{HT}+b$ is a semisimple element of $\g$ with exactly 2
eigenvalues. However in the representation $\Sym^{k-1}$ of $\sl_2$,
every non-zero semisimple element has $k$ distinct eigenvalues.

This completes the proof of Theorem~2, since in that case, $k$ is even
and $V$ has a $\g$-invariant symplectic form, so we must be in case
(iii).

In an earlier version of this paper we gave an ugly proof of the
only case of Proposition 4 needed here ($k$ even, $\g\subset \fraksp(V)$),
involving a detailed case-by-case analysis of minuscule
representations of $\g$.  Subsequently Laumon pointed out to me that
this was a special case of Katz's result, whose proof
also depends on a (longer) case-by-case analysis. I am grateful
to Ian Grojnowski for suggesting a short proof of Katz's general result
along the lines given in the next section.

\section{Lie-theoretic part}

Assume $\g$ satisfies the hypotheses of Proposition~4. 
By the Jacobson-Morosov theorem, there is a
homomorphism $\lambda_{x}\colon \sl_2 \to\g$ such
\[
\lambda_{x}\colon\begin{pmatrix}0&1\\0&0\end{pmatrix}\mapsto x.
\]
Write $\ax\subset\g$ for the image of $\lambda_{x}$, and
\[
y=\lambda_{x}\begin{pmatrix}0&0\\-1&0\end{pmatrix},\quad
h=\lambda_{x}\begin{pmatrix}1&0\\0&-1\end{pmatrix}.
\]
(In the terminology of \cite{Bou}, VIII.11.1, $(x,h,y)$ is an
$\sl_2$-triplet.)

We first observe:

\begin{lem}
$\g$ is simple.
\end{lem}

\begin{proof}
  Suppose that $\g=\g_1\times\g_2$ with $\g_i$ nonzero.  As $V$ is an
  irreducible $\g$-module, it factorises as a tensor product of
  irreducible $\g_i$-modules $V_i$. But since $x$ has maximal rank,
  the restriction of $V$ to $\ax$ is an irreducible representation of
  $\sl_2$, and the tensor product of two non-trivial representations
  of $\sl_2$ is never irreducible, by the Clebsch-Gordan formula.
\end{proof}

The triple $(x,h,y)$ is a principal $\sl_2$-triplet in $\sl(V)$, since
$x$ has maximal rank, and so is also a principal $\sl_2$-triplet in
$\g$ (\cite{Bou} VIII.11.4).  Let $n$ be the rank of $\g$ and $1\le
r_1 \le r_2\le\dots\le r_n$ be the exponents of its root system. Then
one knows (\cite{Ko}, or see for example \cite{Bou} VIII.11, exercise
11) that under the adjoint action of $\ax\simeq\sl_2$, $\g$ decomposes
as the direct sum of the irreducible representations $\Sym^{2r_i}$.

For the adjoint action of $\ax$ on $\sl(V)$ the
exponents are $\{1, 2, \dots,k-1\}$ and one can write down the
decomposition into irreducibles totally explicitly: consider the
matrix powers $x^r\in \sl(V)$ for $1\le r\le k-1$. Let $U_r$
be the $\ax$-submodule of $\sl(V)$ generated by
$x^r$. Obviously $\ad(x)x^r=0$, and since $[h,x]=2x$ one gets
$\ad(h)x^r=2rx^r$. Thus $x^r$ is a highest weight vector in $U_r$,
which is isomorphic to $\Sym^{2r}$, and a basis for $U_r$ is given by
$\{\ad(y)^ix^r\mid 0\le i\le 2r\}$. Thus
\[
\sl(V)= \bigoplus_{r=1}^{k-1} U_r.
\]
Therefore $\g=\bigoplus_{i=1}^n U_{r_i}$. In particular this proves
part (i) of the following Lemma.

\begin{lem}
  (i) The exponents of $\g$ are distinct and satisfy $r_i\le
  k-1$.
  
  (ii) If $r$ and $s$ are exponents of $\g$ and $r+s\le k$ then $r+s-1$
  is also an exponent of $\g$.
\end{lem}

\begin{proof}[of (ii)]
As the Lie bracket $\g\otimes\g \to\g$ is $\ax$-equivariant and $U_r\simeq
\Sym^{2r}$, by the Clebsch-Gordan formula we see that if $r\ge s$ then
\[
[U_r,U_s]=\bigoplus_{t\in T} U_t
\]
for some subset $T\subset \{t\in\Z\mid r-s\le t\le
\min(r+s,k-1)\}$. If $r+s\in T$ then the Lie bracket would give a
non-zero pairing $U_r\otimes U_s\to U_{r+s}$, which would necessarily
be non-zero on the tensor product of the highest weight vectors. But
$[x^r,x^s]=0$, hence $r+s\notin T$. On the other hand, since $x^s$ is
a highest weight vector for $U_s$ one has
$\ad(x)\ad(y)x^s=2sx^s$, and therefore
\[
[x^r,\ad(y)x^s]=x[x^{r-1},\ad(y)x^s]+[x,\ad(y)x^s]x^{r-1}
=x[x^{r-1},\ad(y)x^s]+2sx^{r+s-1}
\]
and so by induction one obtains
\[
[x^r,\ad(y)x^s]=2rs x^{r+s-1}.
\]
Therefore $[U_r,U_s]\supset U_{r+s-1}$ if $r+s\le k$.
As $\g$ is a Lie subalgebra of $\sl(V)$, the Lemma follows. 
\end{proof}

So to finish the proof, it suffices to determine those simple Lie
algebras which admit a representation of dimension $k$ and whose
exponents satisfy the conditions of Lemma~6. From standard tables (for
example \cite{Bou}, Chapters IV and VIII) one extracts the
information contained in the table below.
%\begin{table}
%\caption{}

\begin{center}
\begin{tabular}{lll}
\\
$\g$&exponents of $\g$&least dimensions of representations\\
\hline
$A_n$&1, 2, 3, \dots, $n$&$n+1$, $n(n+1)/2$\\
$B_n$,&1, 3, 5, \dots, $2n-1$&$2n+1$, $n(2n+1)$\\
$C_n$,&1, 3, 5, \dots, $2n-1$&$2n$, $n(2n-1)$\\
$D_n$&1, 3, 5, \dots, $2n-3$, $n-1$\\
$E_6$&1, 4, 5, 7, 8, 11\\
$E_7$&1, 5, 7, 9, 11, 13, 17\\
$E_8$&1, 7, 11, 13, 17, 19, 23, 29\\
$F_4$&1, 5, 7, 11\\
$G_2$&1, 5&7, 14\\
\\
\end{tabular}
\end{center}

%\end{table}
From this one sees that the only cases satisfying the hypotheses of
Lemma~6 are: $A_1$ with $k$ arbitrary; $A_n$ with $k=n+1$; $B_n$ with
$k=2n+1$; $C_n$ with $k=2n$; and $G_2$ with $k=7$, which are precisely
those cases listed in Proposition~4.

%\begin{acknowledgements}
The author would like to thank D. Blasius, G. Harder, and J-P.
Serre for useful discussions, G. Laumon for drawing his attention to
\cite{Ka}, and especially I. Grojnowski for suggesting the proof of
Proposition 4.
%\end{acknowledgements}

\enlargethispage{25mm}\vspace{6mm}\noindent
   Department of Pure Mathematics and Mathematical Statistics\\
   Centre for Mathematical Sciences\\
   Wilberforce Road\\
   Cambridge\enspace CB3 0WB\\
   \email{a.j.scholl@dpmms.cam.ac.uk}%
\end{document}